\documentclass{article}

\begin{document}

\input psbox.tex
\psfordvips

\centerline{\LARGE\bf Simultaneous packing and covering }

\medskip
\centerline{\LARGE\bf in the Euclidean plane II\footnote{This work
is supported by the National Science Foundation of China, a 973
Project and a special grant from Peking University.}}

\bigskip \centerline{\bf Chuanming Zong\footnote{{\it E-mail address:} cmzong@math.pku.edu.cn}}

\bigskip \centerline{\it School of Mathematical Sciences, Peking
University, Beijing 100871, China}

\bigskip
\centerline{\it Dedicated to Professor Edmund Hlawka on the occasion
of his 90th birthday}

\bigskip
\noindent {\bf Abstract}

\medskip
In 1950, C.A. Rogers introduced and studied the simultaneous
packing and covering constants for a convex body and obtained the
first general upper bound. Afterwards, they have attracted the
interests of many authors such as L. Fejes T\'oth, S.S.
Ry$\breve{s}$kov, G.L. Butler, K. B\"or\"oczky, H. Horv\'ath, J.
Linhart and M. Henk since, besides their own geometric significance,
they are closely related to the packing densities and the covering
densities of the convex body, especially to the Minkowski-Hlawka
theorem. However, so far our knowledge about them is still very
limited. In this paper we will determine the optimal upper bound of
the simultaneous packing and covering constants for two-dimensional
centrally symmetric convex domains, and characterize the domains
attaining the upper bound.



\section*{1. Introduction}

\vspace{0.4cm} In 1950, C.A. Rogers introduced and studied two
constants $\gamma (K)$ and $\gamma^*(K)$ for an $n$-dimensional {\it
convex body} $K$. Namely, $\gamma (K)$ is the smallest positive
number $r$ such that there is a translative packing $K+X$ satisfying
$E^n=r K+X$, and $\gamma^* (K)$ is the smallest positive number
$r^*$ such that there is a lattice packing $K+\Lambda$ satisfying
$E^n=r^* K+\Lambda$. In some references, they are called the {\it
simultaneous packing and covering constants} for the convex body.
Clearly, these numbers are closely related to the {\it packing
densities} and the {\it covering densities} of the convex body,
especially to the {\it Minkowski-Hlawka theorem}.

In 1970 and 1978, S.S. Ry$\check{s}$kov  and L. Fejes T\'oth
independently introduced and investigated two related numbers $\rho
(K)$ and $\rho^*(K)$, where $\rho (K)$ is the largest positive
number $r$ such that one can put a translate of $rK$ into every
translative packing $K+X$, and $\rho^* (K)$ is the largest positive
number $r^*$ such that one can put a translate of $r^*K$ into every
lattice packing $K+\Lambda $.

Clearly, for every convex body $K$ we have
$$\gamma (K)\le \gamma^*(K)$$
and
$$\rho (K)\le \rho^*(K).$$
As usual, let $C$ denote an $n$-dimensional centrally symmetric
convex body. Then, we also have
$$\gamma (C)=\rho (C)+1$$
and
$$\gamma^* (C)=\rho^* (C)+1.$$

Let $B^n$ denote the $n$-dimensional unit ball. Just like the {\it
packing density problem} and the {\it covering density problem}, to
determine the values of $\gamma (B^n)$ and $\gamma^*(B^n)$ is
important and interesting. However, so far our knowledge about
$\gamma (B^n)$ and $\gamma^*(B^n)$ is very limited. We list the main
known results in the following table.

$$\begin{tabular}{|c|c|c|c|c|}
\hline
n&2&3&4&5\\
\hline $\gamma^*(B^n)$&$\sqrt{4\over 3}$&$\sqrt{5\over
3}$&$\sqrt{2\sqrt{3}}(\sqrt{3}-1)$&
$\sqrt{{3\over 2}+{\sqrt{13}\over 6}}$\\
\hline Author & &B\"or\"oczky \cite{boro86}&Horv\'ath
\cite{horv82}&Horv\'ath \cite{horv86}\\
\hline
\end{tabular}$$

\medskip Let $\delta (K)$ and $\delta^*(K)$ denote the maximal
{\it translative packing density} and the maximal {\it lattice
packing density} of $K$, respectively. A fundamental problem in
Packing and Covering is to determine {\it if $$\delta
(K)=\delta^*(K)$$ holds for every convex body}. It is easy to see
that $\gamma^*(C)\ge 2$ will imply
$$\delta (C)\ge 2 \delta^*(C),\eqno (1)$$
which will give a negative answer to the previous problem. On the
other hand, if $\gamma^*(C)\le 2-k$ holds for a positive constant
$k$ and for every centrally symmetric convex body $C$, then the
Minkowski-Hlawka theorem can be improved to
$$\delta^*(C)\ge {1\over {(2-k)^n}}.\eqno (2)$$

In 1950, C.A. Rogers \cite{roge50} discovered a constructive method
by which he deduced
$$\gamma^*(C)\le 3$$
for all $n$-dimensional centrally symmetric convex bodies. In 1972,
by mean value techniques developed by C.A. Rogers and C.L. Siegel,
the above upper bound was improved by G.L. Butler \cite{butl72} to
$$\gamma^*(C)\le 2+o(1).$$
This result is fascinating, because it gives hope to both (1) and
(2).

In two and three dimensions, as one can imagine, the situation is
much better. In 1978, based on an ingenious idea of I. F\'ary
\cite{fary50}, J. Linhart \cite{linh78} proved that
$$\gamma (K)=\gamma^*(K)\le \mbox{${3\over 2}$}$$
holds for every two-dimensional convex domain and the upper bound is
attained only by triangles. However, just like the packing density
problem, to determine the best upper bound for $\gamma^* (C)$ turns
out to be much more challenging. Recently C. Zong \cite{zong02} and
\cite{zong03} obtained
$$\gamma^*(C)\le 1.2$$
for all two-dimensional centrally symmetric convex domains and
$$\gamma^*(C)\le 1.75$$
for all three-dimensional centrally symmetric convex bodies.
Needless to say, neither of them is optimal. In this paper we will
prove the following theorem.

\bigskip\noindent {\bf Theorem.} {\it For every two-dimensional
centrally symmetric convex domain $C$ we have
$$\gamma (C)=\gamma ^*(C)\le 2(2-\sqrt{2})\approx 1.17157\cdots ,$$
where the second equality holds if and only if $C$ is an affinely
regular octagon.}

\bigskip\noindent {\bf Remark 1.} {\it To determine the above upper
bound has been listed as an open problem in several references, for
example in Brass, Moser and Pach {\rm \cite{bras05}}, Linhart {\rm
\cite{linh78}}, Zong {\rm \cite{zong02}} and {\rm \cite{zong02'}}.}

\medskip\noindent
{\bf Remark 2.} {\it The identity $$\gamma (C)=\gamma ^*(C)$$ was
proved by J. Linhart {\rm \cite{linh78}} and C. Zong {\rm
\cite{zong02}}. We restate it here just for completion.}

\medskip\noindent
{\bf Remark 3.} {\it Let $\theta (K)$ and $\theta^*(K)$ denote the
least translative covering density and the least lattice covering
density of $K$, respectively. In the plane it was proved by L. Fejes
T\'oth that
$$\theta (C)=\theta^*(C)\le {{2\pi }\over {3\sqrt{3}}},$$
where the second equality holds if and only if $C$ is an ellipse. It
was proved by C.A. Rogers {\rm \cite{roge51}} that $$\delta
(K)=\delta^* (K)$$ holds for every two-dimensional convex domain $K$
in $1951$. However, to find the optimal lower bound for $\delta (C)$
is still a challenging open problem {\rm (}see Reinhardt {\rm
\cite{rein34}}, Mahler {\rm \cite{mahl47}} or Brass, Mosser and Pach
{\rm \cite{bras05})}. Nevertheless, it has been proved that neither
ellipses nor affinely regular octagons can attain the optimal lower
bound.}

\bigskip\medskip\section*{2. Several Basic Lemmas}

\vspace{0.4cm} Let $\partial (K)$ and $int (K)$ denote the {\it
boundary} and the {\it interior} of $K$, respectively. As usual, we
call a convex body {\it regular} if for every point ${\bf x}\in
\partial (K)$ there is a unique tangent hyperplane and every tangent
plane touches its boundary at a single point. For convenience, in
the rest of this paper $C$ always means a two-dimensional centrally
symmetric convex domain. Now let us introduce several basic lemmas
which will be useful in our proof.

\bigskip\noindent
{\bf  Lemma 1 (Mahler \cite{mahl47}).} {\it If $\pm {\bf v}_1$, $\pm
{\bf v}_2$ and $\pm {\bf v}_3$ are the six vertices of an affinely
regular hexagon inscribed in $C$, then $C+\Lambda $ is a lattice
packing of $C$, where $$\Lambda =\{ 2z_1{\bf v}_1+2z_2{\bf v}_2:\
z_i\in Z\}.$$}

\medskip\noindent
{\bf Lemma 2 (Eggleston \cite{eggl58}).} {\it For every convex body
there is a sequence of regular convex bodies which converges to the
convex body in the sense of Hausdorff metric.}

\bigskip Let ${\bf v}_1, {\bf v}_2, \cdots , {\bf v}_6 $ be the six vertices (in
anti-clock order) of a centrally symmetric hexagon $H$ which is
inscribed in $C$, let ${\bf m}_i$ denote the midpoint of ${\bf
v}_i{\bf v}_{i+1}$, and let ${\bf m}^*_i$ denote the point in the
direction of ${\bf m}_i$ and on the boundary of $C$. Then, for $i=1,
2$ and $3$,  we define
$$f_i({\bf v}_1)={{\|{\bf o}, {\bf m}^*_i\|}\over
{\| {\bf o}, {\bf m}_i\| }},$$ where $\| {\bf x}, {\bf y}\|$ denotes
the {\it Euclidean distance} between ${\bf x}$ and ${\bf y}$.

\bigskip
\noindent {\bf Lemma 3 (Zong \cite{zong96}).} {\it For any ${\bf
x}\in \partial (C)$, we can choose five points ${\bf x}_2,$  ${\bf
x}_3$, ${\bf x}_4$, ${\bf x}_5$ and ${\bf x}_6$ from $\partial (C)$
such that they together with ${\bf x}$ are the six vertices of an
affinely regular hexagon. When $C$ is regular and ${\bf x}$ moves
along $\partial (C)$, we can choose the points such that all
$f_1({\bf x})$, $f_2({\bf x})$ and $f_3({\bf x})$ are continuous
functions of ${\bf x}$.}

\medskip\noindent
{\bf Proof.} In fact, \cite{zong96} only contains a proof for the
first part of this lemma. Here let us outline a proof for the second
part.

Assume that $C$ is regular and, without loss of generality, we
assume further that ${\bf x}{\bf x}_2{\bf x}_3{\bf x}_4{\bf x}_5{\bf
x}_6$ is a regular hexagon with ${\bf x}=(1,0)$, as shown in Figure
1. Let $\epsilon $ be a positive number, let ${\bf x}'$ be a point
on the boundary of $C$ such that $\angle {\bf x}'{\bf o} {\bf
x}=\epsilon ,$ let $\Gamma_2$ and $\Gamma_3$ denote the straight
lines which are parallel with ${\bf o}{\bf x}'$ and pass ${\bf x}_2$
and ${\bf x}_3$, respectively.

$$\psannotate{\psboxto(0cm;7cm){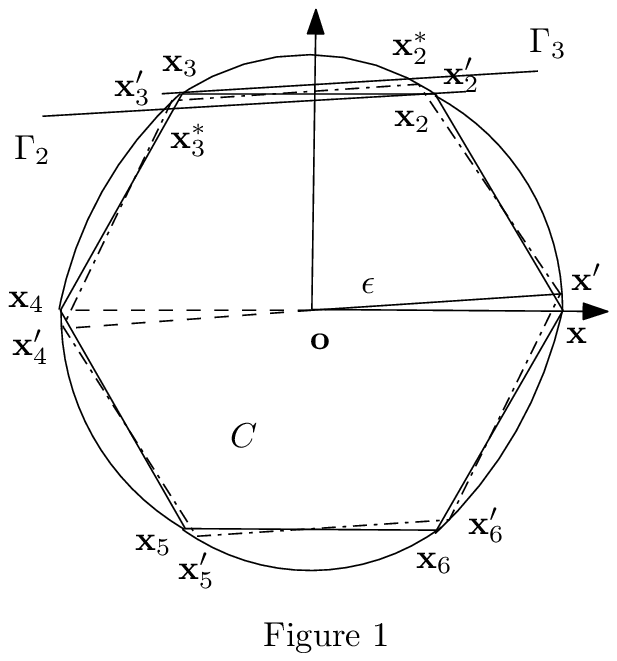}}{}$$

Clearly, when $\epsilon $ is sufficiently small, $\Gamma_2$
intersects $\partial (C)$ at two points ${\bf x}_2$ and ${\bf
x}^*_3$, $\Gamma_3$ intersects $\partial (C)$ at two points ${\bf
x}_3$ and ${\bf x}^*_2$, and both $\| {\bf x}_2, {\bf x}^*_2\|$ and
$\| {\bf x}_3, {\bf x}^*_3\|$ are small. In addition, then the three
directions ${\bf x}{\bf x}'$, ${\bf x}_2{\bf x}^*_2$ and ${\bf
x}_3{\bf x}^*_3$ are approximately the tangent directions of $C$ at
${\bf x}$, ${\bf x}_2$ and ${\bf x}_3$, respectively. Thus,
comparing triangles ${\bf o}{\bf x}{\bf x}'$ and ${\bf x}_2{\bf
x}_3{\bf x}_3^*$ with ${\bf x}_3{\bf x}_2{\bf x}_2^*$ and ${\bf
o}{\bf x}_4{\bf x}'_4$ (where ${\bf x}'_4=-{\bf x}'$), respectively,
by convexity and elementary geometry we get
$$\| {\bf x}^*_2, {\bf x}_3\|< \| {\bf o}, {\bf x}'\| <\| {\bf x}_2, {\bf x}^*_3\|$$
when $\epsilon $ is sufficiently small.
Therefore $\partial (C)$ has two points ${\bf x}'_2$ and ${\bf
x}'_3$ between $\Gamma_2$ and $\Gamma_3$ such that ${\bf x}'_3{\bf
x}'_2$ is parallel with ${\bf o}{\bf x}'$ and
$$\| {\bf x}'_3, {\bf x}'_2\|=\| {\bf o}, {\bf x}'\| .$$

Taking ${\bf x}'_4=-{\bf x}'$, ${\bf x}'_5=-{\bf x}'_2$ and ${\bf
x}'_6=-{\bf x}'_3$, it is easy to see that ${\bf x}'{\bf x}'_2{\bf
x}'_3{\bf x}'_4{\bf x}'_5{\bf x}'_6$ is an affinely regular hexagon.
Since both ${\bf x}'_2$ and ${\bf x}'_3$ continuously depend on
${\bf x}'$, all $f_1({\bf x})$, $f_2({\bf x})$ and $f_3({\bf x})$
are continuous functions of ${\bf x}$. The lemma is proved.

\medskip\noindent
{\bf Remark 4.} {\it Without the regularity assumption the second
part of the lemma will not be true. For example, when
$$C=\{ (x,y):\ |x|\le 1,\ |y|\le 1\},$$ the corresponding function $f_1({\bf x})$
is not continuous at ${\bf x}=(1,0)$.}

\medskip\noindent
{\bf Lemma 4 (Zong \cite{zong02}).} {\it Let $f_1({\bf v}_1)$,
$f_2({\bf v}_1)$ and $f_3({\bf v}_1)$ be the numbers defined above
Lemma $3$. Then we have}
$$\gamma^*(C)\le \max \{ f_1({\bf v}_1), f_2({\bf v}_1), f_3({\bf v}_1)\} .$$

\bigskip\medskip
\section*{ 3. A Proof for the Theorem}

\vspace{0.4cm} For convenience, let $L({\bf x}, {\bf y})$ denote the
straight line passing two points ${\bf x}$ and ${\bf y}$, and write
$$\alpha =2(2-\sqrt{2}).$$ To make the complicated proof more
transparent, we divide it into three parts.

\bigskip\noindent
{\bf Assertion I.} {\it For every two-dimensional centrally
symmetric convex domain $C$, there is a corresponding inscribed
affinely regular hexagon ${\bf v}_1{\bf v}_2\cdots {\bf v}_6$
satisfying}
$$f_1({\bf v}_1)\ge f_2({\bf v}_1)=f_3({\bf v}_1).$$

\smallskip\noindent
{\bf Proof.} First let us consider the case that $C$ is regular. By
Lemma 3, all $f_1({\bf x})$, $f_2({\bf x})$ and $f_3({\bf x})$ are
continuous functions of ${\bf x}\in \partial (C)$. Therefore,
$$f({\bf x})= \min_{i=1,2}\{ f_i({\bf x})\} -f_3({\bf x})$$ is
also a continuous function of ${\bf x}\in \partial (C)$. If, without
loss of generality,
$$f_i({\bf x}_1)> f_3({\bf x}_1), \quad i=1,2,$$
hold at some point ${\bf x}_1\in \partial (C)$ and ${\bf x}_1{\bf
x}_2{\bf x}_3{\bf x}_4{\bf x}_5{\bf x}_6$ is the corresponding
affinely regular hexagon inscribed in $C$, then we get
$$f({\bf x}_1)=\min_{i=1,2}\{ f_i({\bf x}_1)\} -f_3({\bf x}_1)>0$$ and
$$f({\bf x}_2)=f_3({\bf x}_1)-f_1({\bf x}_1)<0.$$
Therefore, there are two suitable points ${\bf v}, {\bf v}_1\in
\partial (C)$ satisfying
$$f({\bf v})= \min_{i=1,2}\{ f_i({\bf v})\} -f_3({\bf v})=0$$
and
$$f_1({\bf v}_1)\ge f_2({\bf v}_1)=f_3({\bf v}_1).$$

In the general case, by Lemma 2, there is a sequence of regular
centrally symmetric convex domains $C_1$, $C_2$, $\cdots $ which
converges to $C$ in the sense of the Hausdorff-metric. Assume that
$$H_i={\bf v}^i_1{\bf v}^i_2{\bf v}^i_3{\bf v}^i_4{\bf v}^i_5{\bf
v}^i_6$$ is an affinely regular hexagon inscribed in $C_i$ and
satisfying
$$f_1({\bf v}^i_1)\ge f_2({\bf v}^i_1)=f_3({\bf v}^i_1).$$
Then, by {\it Blaschke's selection theorem}, there is a subsequence
of the sequence $H_1$, $H_2$, $\cdots $ which converges to an
affinely regular hexagon ${\bf v}_1{\bf v}_2{\bf v}_3{\bf v}_4{\bf
v}_5{\bf v}_6$ which is inscribed in $C$ and satisfies
$$f_1({\bf v}_1)\ge f_2({\bf v}_1)=f_3({\bf v}_1).$$
Thus, Assertion I is proved.

\bigskip\noindent
{\bf Assertion II.} {\it For each two-dimensional centrally
symmetric convex domain $C$ there is a corresponding lattice
$\Lambda$ such that $C+\Lambda$ is a packing and $\alpha C+\Lambda$
is a covering in $E^2$.}

\medskip\noindent
{\bf Proof.} Let ${\bf v}_1{\bf v}_2\cdots {\bf v}_6$ be the hexagon
obtained in Assertion I. For convenience, we write
$$\kappa =f_1({\bf v}_1),\quad \lambda =f_2({\bf v}_1)$$ and define
$$\Lambda_1 =\{ 2z_2{\bf v}_2+2z_3{\bf v}_3:\ z_i\in Z\}.$$
By Lemma 1, it follows that $C+\Lambda_1 $ is a lattice packing.

If $\kappa < \alpha $, then by Lemma 4 we can get
$$\gamma^*(C)\le \kappa < \alpha .$$
Thus, from now on we assume that $\kappa \ge \alpha $.

$$\psannotate{\psboxto(0cm;8cm){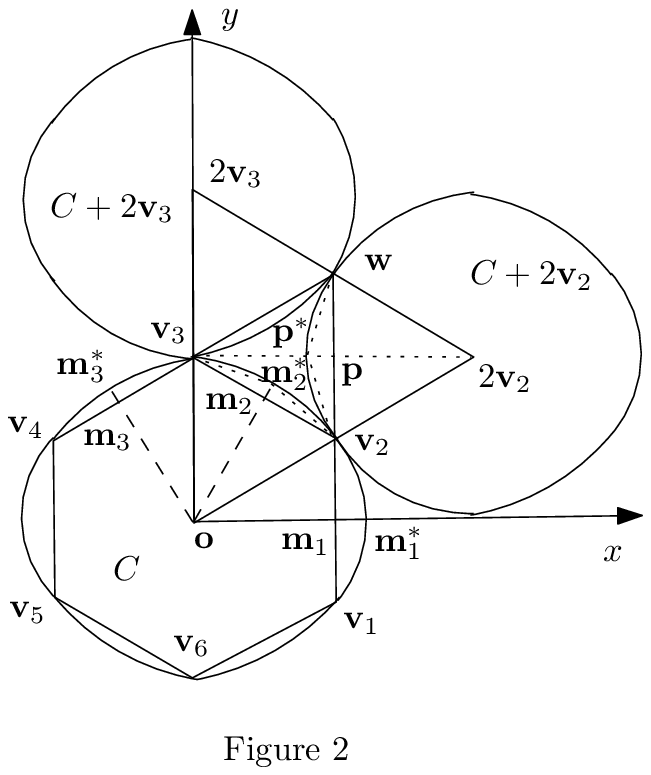}}{}$$

As it is shown in Figure 2, without loss of generality, we assume
that ${\bf v}_1{\bf v}_2{\bf v}_3{\bf v}_4{\bf v}_5{\bf v}_6$ is a
regular hexagon with ${\bf v}_2=(\sqrt{3}/2, 1/2)$ and ${\bf
v}_3=(0,1)$.  By a routine computation based on elementary geometry
it can be shown that the equation of $L(\alpha {\bf v}_3, \alpha
{\bf m}^*_2)$ is
$$y-\alpha ={{3\lambda -4}\over {\sqrt{3}\lambda }}x$$
and therefore
$$L(\alpha {\bf v}_3, \alpha {\bf m}^*_2)\cap L({\bf v}_3, 2{\bf
v}_2)=\left( {{\sqrt{3}\lambda (\alpha -1)}\over {4-3\lambda }},
1\right).\eqno (3)$$ Let ${\bf p}$ denote the midpoint of ${\bf
v}_2{\bf w}$ and let ${\bf p}^*$ denote the boundary point of
$C+2{\bf v}_2$ in the direction from $2{\bf v}_2$ to ${\bf p}$. By
symmetry we have
$${{\| 2{\bf v}_2, {\bf p}^*\|}\over {\| 2{\bf v}_2, {\bf
p}\|}}={{\| {\bf o}, {\bf m}^*_1\|}\over {\| {\bf o}, {\bf
m}_1\|}}=\kappa .$$ Therefore, by a routine computation it can be
deduced that
$$[{\bf v}_3, 2{\bf v}_2] \cap \partial (\alpha C+2{\bf v}_2)=\left(
\sqrt{3}\left(1-\mbox{${1\over 2} $} \alpha \kappa
\right),1\right),\eqno (4)$$ where $[{\bf x}, {\bf y}]$ denotes the
segment between ${\bf x}$ and ${\bf y}$. Then, by (3) and (4) it
follows that $\alpha C+\Lambda_1 $ will be a covering of $E^2$ and
therefore
$$\gamma ^*(C)< \alpha $$ if
$$\sqrt{3}\left(1-{{\alpha \kappa }\over 2}\right)< {{\sqrt{3}\lambda (\alpha
-1)}\over {4-3\lambda }}.\eqno (5)$$

By convexity it is easy to see from Figure 2 that $\lambda \le 4/3$.
Then it can be easily verified that (5) holds whenever $\kappa
>\sqrt{2}$ or $\kappa \ge \lambda >\alpha $. Thus, in the rest of the proof we assume that
$$1\le \lambda \le \alpha, \eqno (6)$$

$$\alpha \le \kappa \le \sqrt{2}\eqno (7)$$
and
$${{4-3\lambda }\over \lambda }\ge {{\alpha -1}\over {1-{1\over 2}\alpha \kappa }}.\eqno (8)$$

$$\psannotate{\psboxto(0cm;7.8cm){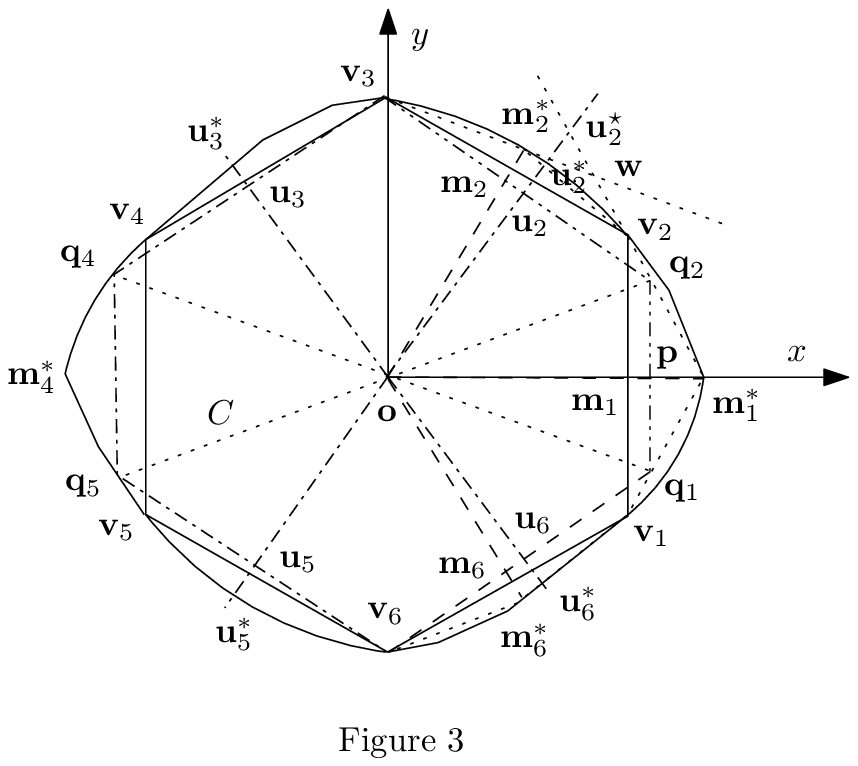}}{}$$

As it is shown in Figure 3, without loss of generality, we assume
that ${\bf v}_1{\bf v}_2{\bf v}_3{\bf v}_4{\bf v}_5{\bf v}_6$ is a
regular hexagon with ${\bf v}_2=(\sqrt{3}/2, 1/2)$ and ${\bf
v}_3=(0,1)$. Since $\| {\bf o}, {\bf m}^*_1\| / \| {\bf o}, {\bf
m}_1\| =\kappa $ and $\| {\bf o}, {\bf m}^*_2\| / \| {\bf o}, {\bf
m}_2\| =\lambda $, by a routine computation we get that the
equations of $L({\bf v}_3, {\bf m}_2^*)$ and $L({\bf v}_2, {\bf
m}_1^*)$ are
$$y-1={{3\lambda -4}\over {\sqrt{3}\lambda }}x\eqno (9)$$
and
$$y={{x-{\sqrt{3}\over 2}\kappa }\over {\sqrt{3}(1-\kappa )}},\eqno (10)$$
respectively.

Since $\| {\bf o}, {\bf m}^*_1\| /\| {\bf o}, {\bf m}_1\| \ge \alpha
$, there is a point ${\bf p}$ between ${\bf m}_1$ and ${\bf m}_1^*$
satisfying $${{\| {\bf o}, {\bf m}^*_1\|}\over {\| {\bf o}, {\bf
p}\| }}=\alpha .\eqno (11)$$ Then, there are two points ${\bf
q}_1\in L({\bf v}_1, {\bf m}_1^*)$ and ${\bf q}_2\in L({\bf v}_2,
{\bf m}_1^*)$ such that ${\bf p}\in L({\bf q}_1, {\bf q}_2)$ and
$L({\bf q}_1, {\bf q}_2)$ is parallel with $L({\bf v}_1, {\bf
v}_2)$. It follows by a routine computation that
$${\bf p}=\left( {{\sqrt{3}\kappa }\over {2\alpha }}, 0\right)$$
and $${\bf q}_2=\left( {{\sqrt{3}\kappa }\over {2\alpha }}, {{\kappa
(\alpha -1)}\over {2(\kappa -1)\alpha }}\right).$$

Let ${\bf u}_2$ denote the midpoint of ${\bf v}_3{\bf q}_2$. It is
easy to see that
$${\bf u}_2=\left( {{\sqrt{3}\kappa }\over {4\alpha }}, {{\kappa (\alpha
-1)}\over {4(\kappa -1)\alpha }}+{1\over 2}\right)$$ and the
equation of $L({\bf o}, {\bf u}_2)$ is
$$y=\left( {{\kappa (\alpha -1)}\over {4(\kappa -1)\alpha }}+{1\over 2}\right)
{{4\alpha }\over {\sqrt{3}\kappa }} x.\eqno (12)$$ Let ${\bf u}_2^*$
and ${\bf u}_2^\star$ denote the intersections of $L({\bf o}, {\bf
u}_2)$ with $L({\bf v}_3, {\bf m}_2^*)$ and $L({\bf v}_2, {\bf
m}_1^*)$, respectively. By (9), (10) and (12), we can get
$${\bf u}_2^*=\left( \left({{\alpha -1}\over
{\sqrt{3}(\kappa -1)}}+{{2\alpha }\over {\sqrt{3}\kappa
}}+{{4-3\lambda }\over {\sqrt{3}\lambda }}\right)^{-1}, y^*\right)$$
and
$${\bf u}_2^\star =\left( {{\sqrt{3}\kappa ^2}\over {2\alpha (3\kappa -2)}},
y^\star \right),$$ where the $y$-coordinates of both ${\bf u}_2^*$
and ${\bf u}_2^\star$ are not necessary for our purpose. Thus, we
get
$${{\| {\bf o}, {\bf u}_2^*\|}\over {\| {\bf o}, {\bf
u}_2\|}}={{4\alpha }\over \kappa }\left( {{\alpha -1}\over {\kappa
-1}}+{{2\alpha }\over \kappa }+{{4-3\lambda }\over \lambda
}\right)^{-1}$$ and
$${{\| {\bf o}, {\bf u}_2^\star\|}\over {\| {\bf o}, {\bf
u}_2\|}}={{2\kappa }\over {3\kappa -2}}.$$  For convenience, we
write
$$f(\kappa ,\lambda )={{4\alpha }\over \kappa }\left( {{\alpha -1}\over
{\kappa -1}}+{{2\alpha }\over \kappa }+{{4-3\lambda }\over \lambda
}\right)^{-1}$$ and
$$g(\kappa )={{2\kappa }\over {3\kappa -2}}.$$

Next, we proceed to show $f(\kappa ,\lambda )\le \alpha$ and
$g(\kappa )\ge \alpha $. It follows by (8) that
$$f(\kappa ,\lambda )\le {{4\alpha }\over \kappa }\left( {{\alpha -1}\over
{\kappa -1}}+{{2\alpha }\over \kappa }+{{\alpha -1}\over {1-{1\over
2}\alpha \kappa }}\right)^{-1}.$$ By routine computations, it is
easy to see that
$${{4\alpha }\over \kappa }\left( {{\alpha -1}\over {\kappa -1}}+{{2\alpha }\over
\kappa }+{{\alpha -1}\over {1-{1\over 2}\alpha \kappa
}}\right)^{-1}\le \alpha
$$ is equivalent with
$${{2\alpha-4}\over \kappa }+{{\alpha -1}\over {\kappa -1}}+{{\alpha
-1}\over {1-{1\over 2}\alpha \kappa }}\ge 0$$ and, substituting
$\alpha $ by $2(2-\sqrt{2})$,
$${{(\sqrt{2}-1)(3-\sqrt{2})^2\left(\kappa -{2\over {3-\sqrt{2}}}\right)^2}\over
{\kappa (\kappa -1)(1-(2-\sqrt{2})\kappa )}}\ge 0,$$ which is
clearly true under the assumption of (6) and (7). Thus, we get
$$f(\kappa ,\lambda )={{\| {\bf o}, {\bf u}_2^*\|}\over {\| {\bf o}, {\bf
u}_2\|}}\le \alpha, \eqno (13)$$ where the equality holds if and
only if $$\left\{
\begin{array}{ll}
\kappa = {2\over {3- \sqrt{2}}}\approx 1.261203875\cdots  &\\
&\\ \lambda = {{2+4\sqrt{2}}\over 7}\approx 1.09383632\cdots . &
\end{array}\right. \eqno(14) $$
On the other hand, it is easy to see that $g(\kappa )$ is a
decreasing function of $\kappa $ when $\kappa $ satisfies (7). Thus,
we have
$$g(\kappa )={{\| {\bf o}, {\bf u}_2^\star\|}\over {\| {\bf o}, {\bf
u}_2\|}}\ge {2\over {3-\sqrt{2}}}>\alpha. \eqno (15)$$

$$\psannotate{\psboxto(0cm;8cm){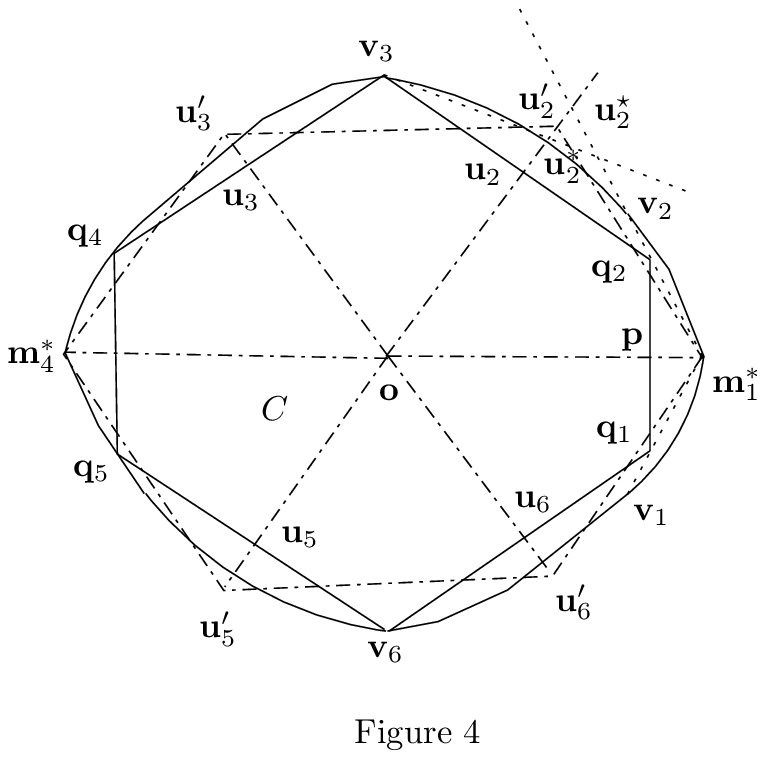}}{}$$

As it is shown in Figure 4, let ${\bf u}'_2$, ${\bf u}'_3$, ${\bf
u}'_5$ and ${\bf u}'_6$ be the four points satisfying
$${{\| {\bf o}, {\bf u}'_i\|}\over {\| {\bf o}, {\bf
u}_i\|}}=\alpha, \quad i=2,\ 3,\ 5,\ 6.\eqno (16)$$ Then ${\bf
m}^*_1{\bf u}'_2{\bf u}'_3{\bf m}^*_4{\bf u}'_5{\bf u}'_6$ is an
affinely regular hexagon. For convenience, we write
$$H=conv\{ {\bf m}^*_1, {\bf u}'_2, {\bf u}'_3, {\bf m}^*_4, {\bf u}'_5, {\bf
u}'_6\},$$

$$C'=conv\{ C, {\bf u}'_2, {\bf u}'_3, {\bf u}'_5, {\bf u}'_6\}$$
and $$\Lambda_2 =\{ 2z_1{\bf m}^*_1+2z_2{\bf u}'_2:\ z_i\in Z\}.$$
By (13), (15) and convexity it follows that
$$\{ {\bf m}^*_1, {\bf u}'_2, {\bf u}'_3, {\bf m}^*_4, {\bf u}'_5, {\bf
u}'_6\}\subset \partial (C').$$ Thus, by Lemma 1, $C'+\Lambda_2 $ is
a packing and therefore $C+\Lambda_2 $ is a packing too. On the
other hand, it follows by (16) that $\alpha H+\Lambda_2 $ is a
tiling of $E^2$ and therefore $\alpha C+\Lambda_2 $ is a covering of
$E^2$. Thus, we get
$$\gamma ^*(C)\le \alpha =2(2-\sqrt{2}).\eqno (17)$$
Assertion II is proved.

\bigskip\noindent
{\bf Assertion III.} {\it The equality
$$\gamma^*(C)=2(2-\sqrt{2})$$ holds if and only if $C$ is an affinely regular octagon.}

\medskip\noindent
{\bf Proof.} Let $P_8$ denote an affinely regular octagon. It was
proved by Linhart \cite{linh78} and Zong \cite{zong02} that
$$\gamma^* (P_8)=2(2-\sqrt{2}).$$
On the other hand, if $D$ is a two-dimensional centrally symmetric
convex domain satisfying
$$\gamma ^*(D)=2(2-\sqrt{2}),\eqno (18)$$
we proceed to show that it must be an affinely regular octagon.

\medskip
First of all, by reexamining the proof of Assertion II, especially
(13) and the construction to prove (17), it is not hard to see that
(18) holds only if the corresponding $\kappa $ and $\lambda $
satisfy (14).

\medskip
Second, using the notation in Figure 3, we claim that
$${\bf v}_2{\bf m}_1^*\subset \partial (D)\eqno (19)$$
and
$${\bf v}_3{\bf u}_2^*\subset \partial (D).\eqno (20)$$
If, on the contrary, (19) does not hold, then
$$(1+\epsilon ){\bf q}_2\in int (D)$$
holds with small positive number $\epsilon $. For convenience, we
write $${\bf q}'_2=(1+\epsilon ){\bf q}_2.$$ Then there is a point
${\bf q}'_1\in ({\bf q}_1, {\bf v}_1)$ such that the midpoint ${\bf
p}'$ of ${\bf q}'_1{\bf q}'_2$ is in $({\bf p}, {\bf m}_1^*)$, where
$({\bf x}, {\bf y})$ denotes the open segment between ${\bf x}$ and
${\bf y}$.  Let ${\bf b}_2$ denote the midpoint of ${\bf v}_3{\bf
q}'_2$, let ${\bf b}_6$ denote the midpoint of ${\bf v}_6{\bf q}'_1$
and let ${\bf b}_i^*$ denote the point on the boundary of $D$ and in
the direction of ${\bf b}_i$. By elementary geometry and convexity
one can deduce that
$${{\| {\bf o}, {\bf m}^*_1\|}\over {\| {\bf o}, {\bf p}'\|}}<
{{\| {\bf o}, {\bf m}^*_1\|}\over {\| {\bf o}, {\bf p}\|}} =\alpha
,$$
$${{\| {\bf o}, {\bf b}^*_2\|}\over {\| {\bf o}, {\bf b}_2\|}}<
{{\| {\bf o}, {\bf u}^*_2\|}\over {\| {\bf o}, {\bf u}_2\|}} =\alpha
\eqno (21)$$ and
$${{\| {\bf o}, {\bf b}^*_6\|}\over {\| {\bf o}, {\bf b}_6\|}}<
{{\| {\bf o}, {\bf u}^*_6\|}\over {\| {\bf o}, {\bf u}_6\|}} =\alpha
.$$ For example, (21) can be deduced from the elementary geometry
illustrated by Figure 5.
$$\psannotate{\psboxto(0cm;6cm){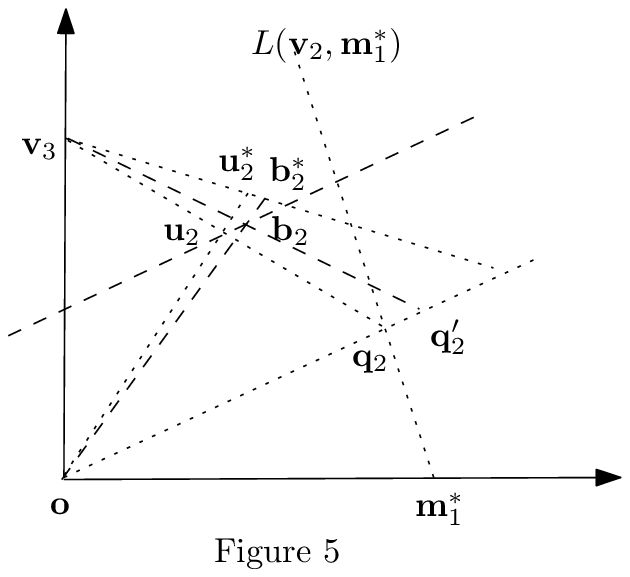}}{}$$

\noindent Then, by a construction similar to that in the proof of
(17) we get
$$\gamma ^*(D)<2(2-\sqrt{2}),$$
which contradicts (18). Thus (19) is proved. The relation (20) can
be shown in a similar way.

\medskip
Finally, let us complete the proof of the assertion based on the
next figure.

$$\psannotate{\psboxto(0cm;6.5cm){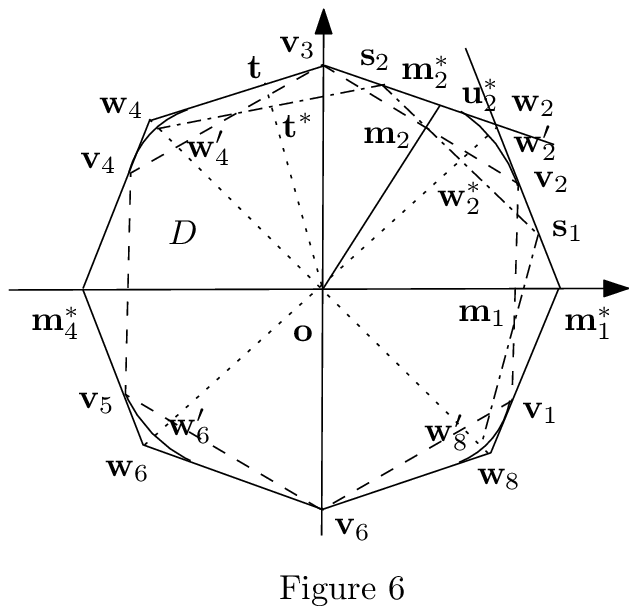}}{}$$

Let ${\bf w}_2$ denote the intersection of $L({\bf v}_3,{\bf
m}^*_2)$ and $L({\bf v}_2, {\bf m}_1^*),$ and let ${\bf w}_4$, ${\bf
w}_6$ and ${\bf w}_8$ be the points defined similarly as shown in
Figure 6. It is easy to verify that ${\bf m}_1^*{\bf w}_2{\bf
v}_3{\bf w}_4{\bf m}_4^*{\bf w}_6{\bf v}_6{\bf w}_8$ is an affinely
regular octagon. In addition, in this case we have ${\bf
v}_3=(0,1)$, ${\bf v}_2=({\sqrt{3}\over 2}, {1\over 2})$, ${\bf
m}_1^*=\left( {{\sqrt{3}}\over {3-\sqrt{2}}}, 0\right),$ ${\bf
u}_2^*=\left( {{\sqrt{3}}\over {2(3-\sqrt{2})}}, {{3-\sqrt{2}}\over
2}\right)$ and ${\bf w}_2=\left( {{\sqrt{3}}\over {3\sqrt{2}-2}},
{\sqrt{2}\over 2}\right).$ We proceed to show
$$D={\bf m}_1^*{\bf w}_2{\bf v}_3{\bf w}_4{\bf m}_4^*{\bf w}_6{\bf
v}_6{\bf w}_8.\eqno (22)$$

If $D$ is not the octagon and let ${\bf w}'_2$ and ${\bf w}'_4$
denote the intersections of $\partial (D)$ with ${\bf o}{\bf w}_2$
and ${\bf o}{\bf w}_4$, respectively. In addition, we write
$$\rho =\| {\bf o}, {\bf w}_2\|/\| {\bf o}, {\bf w}'_2\|,$$
$$\sigma =\| {\bf o}, {\bf w}_4\|/\| {\bf o}, {\bf w}'_4\|$$
and assume that $\rho \ge \sigma $. Based on the coordinates of
${\bf v}_2$, ${\bf u}_2^*$ and ${\bf w}_2$, by a routine computation
we get
$$1\le \sigma\le \rho \le {{4\sqrt{2}+2}\over 7}.\eqno (23) $$

Let ${\bf w}_2^*$ be the point defined by ${\bf o}{\bf w}'_2/{\bf
o}{\bf w}^*_2=\alpha $, let $\Gamma $ denote the straight line
passing ${\bf w}^*_2$ and parallel with $L({\bf w}_4, {\bf w}_8)$,
let ${\bf s}_1$ and ${\bf s}_2$ denote the intersections of $\Gamma
$ with $L({\bf v}_2,{\bf m}_1^*)$ and $L({\bf v}_3, {\bf m}^*_2)$,
respectively, let ${\bf t}^*$ denote the midpoint of ${\bf w}'_4{\bf
s}_2$, let ${\bf t}$ denote the intersection of $L({\bf o}, {\bf
t}^*)$ with $L({\bf v}_3, {\bf w}_4)$, and finally define
$$h(\rho, \sigma)={{\| {\bf o}, {\bf t}\| }\over {\| {\bf o}, {\bf
t}^*\|}}.$$

Based on (23), it can be verified by routine arguments and
computations that
$${\bf s}_1\in {\bf v}_2{\bf m}^*_1\subset \partial
(D),$$ $${\bf s}_2\in {\bf v}_3{\bf u}^*_2\subset \partial (D),$$
$$\| {\bf w}^*_2, {\bf s}_1\|= \| {\bf w}^*_2, {\bf s}_2\|$$ and ${\bf
t}\in \partial (D)$. Especially, by routine but complicated
computations we get
$$h(\rho , \sigma)=\left( {{\sqrt{2}+1}\over 2}-{{2+\sqrt{2}}\over
{4\rho }}+{1\over {2\sigma }}\right)^{-1}.$$ Then, by (23) we get
\begin{eqnarray*}
h(\rho , \sigma)&\le &\left( {{\sqrt{2}+1}\over
2}-{{2+\sqrt{2}}\over
{4\rho }}+{1\over {2\rho }}\right)^{-1}\\
&= & \left( {{\sqrt{2}+1}\over 2}-{{\sqrt{2}}\over
{4\rho}}\right)^{-1}\\
&\le & 2(2-\sqrt{2}),
\end{eqnarray*}
where the final equality holds if and only if $\rho =\sigma =1.$
Then, (22) follows from Lemma 4. Assertion III is proved.

\medskip
As a conclusion of Assertion II and Assertion III the theorem is
proved.

\bigskip\medskip
\section*{ 4. Three Further Remarks}

\vspace{0.4cm}\noindent {\bf Remark 5.} {\it Let
$\lambda_i(C,\Lambda )$ denote the $i$-th {\it successive minimum}
of $C$ with respect to a lattice $\Lambda $, and let $\mu_i
(C,\Lambda )$ denote the $i$-th {\it covering minimum} of $C$ with
respect to $\Lambda$ (see Gruber and Lekkerkerker {\rm
\cite{grub87}} and Kannan and Lov\'asz {\rm \cite{kann88}},
respectively). As a corollary of the theorem we get
$$\min_{\Lambda }{{\mu_2 (C,\Lambda )}\over {\lambda_1(C,\Lambda
)}}\le 2(2-\sqrt{2}),$$ where the equality holds if and only if $C$
is an affinely regular octagon.}

\bigskip\noindent
{\bf Remark 6.} {\it It is well known {\rm (}see L. Fejes T\'oth
{\rm \cite{feje72})} that $${{\theta^* (C)}\over {\delta^* (C)}}\le
{4\over 3}\approx 1.33333\cdots $$ holds for every two-dimensional
centrally symmetric convex domain. However, although our theorem is
optimal, it only can produce
$$ {{\theta^* (C)}\over {\delta^* (C)}}\le \min_{\Lambda }\left(
{{\mu_2 (C,\Lambda )}\over {\lambda_1(C,\Lambda )}}\right)^2\le
8(3-2\sqrt{2})\approx 1.37258\cdots .$$ The reason for this
phenomenon is the optimal covering lattice of a regular octagon is
not homothetic to its optimal packing lattice.}

\bigskip\noindent
{\bf Remark 7.} {\it Let $m_2(C)$ denote the Steiner ratio of the
Minkowski plane determined by a two-dimensional centrally symmetric
convex domain $C$. It is known {\rm (}see Cieslik {\rm
\cite{cies01})} that
$$m_2(C)\le {3\over 4}\ \gamma^*(C).$$
Thus, we have
$$m_2(C)\le {3\over {2+\sqrt{2}}}.$$}

\vspace{1cm} \noindent {\Large\bf Acknowledgements.} In 1996, I
learned this problem from Professor C.A. Rogers when I was a visitor
at University College London. In 2003, when I published two papers
(\cite{zong02} and \cite{zong03}) on this problem, I received an
offprint of \cite{linh78} from Professor J. Linhart. Clearly, he was
not aware of Rogers and Butler's papers on this topic when he
published \cite{linh78}, just like my unawareness of his paper.
Fortunately, our papers have almost no important overlap. I am very
grateful to Professor Rogers for driving my attention to this
problem and to Professor Linhart for sending me his related papers.
For some helpful comments on this paper, I am obliged to Professor
Martin Henk.

\bigskip\medskip
\bibliographystyle{amsplain}

\begin{thebibliography}{99}
\bibitem{boro86}K. B\"or\"oczky, Closest packing and loosest
covering of the space with balls, {\it Studia Sci. Math. Hungar.}
{\bf 21} (1986), 79-89.
\bibitem{bras05}P. Brass, W. Moser and J. Pach, {\it Research
Problems in Discrete Geometry}, Springer-Verlag, New York, 2005.
\bibitem{butl72}G.L. Butler, Simultaneous packing and covering in
Euclidean space, {\it Proc. London Math. Soc.} {\bf 25} (1972),
721-735.
\bibitem{cies01}D. Cieslik, {\it The Steiner Ratio}, Kluwer Academic
Publishers, Dordrecht, 2001.
\bibitem{eggl58}H.G. Eggleston, {\it Convexity}, Cambridge
University Press, Cambridge, 1958.
\bibitem{fary50}I. F\'ary, Sur la densit\'e des r\'eseaux de
domaines convexes, {\it Bull. Soc. Math. France} {\bf 78} (1950),
152-161.
\bibitem{feje93}G. Fejes T\'oth and W. Kuperberg, Packing and
covering with convex sets, {\it Handbook of Convex Geometry} (eds.
P.M. Gruber and J.M. Wills), North-Holland, 1993, 799-860.
\bibitem{feje72}L. Fejes T\'oth, {\it Lagerungen in der Ebene auf
der Kugel und im Raum}, Springer-Verlag, Berlin, 1972.
\bibitem{feje78}L. Fejes T\'oth, Remarks on the closest packing of
convex discs, {\it Comment. Math. Helv.} {\bf 53} (1978), 536-541.
\bibitem{grub87}P.M. Gruber and C.G. Lekkerkerker, {\it Geometry of
Numbers}, North-Holland, Amsterdam, 1987.
\bibitem{henk95}M. Henk, {\it Finite and Infinite Packings},
Habilitationsschrift, Universit\"at Siegen, 1995.
\bibitem{horv82}J. Horv\'ath, On close lattice packing of unit
spheres in the space $E^n$, {\it Proc. Steklov Math. Inst.} {\bf
152} (1982), 237-254.
\bibitem{horv86}J. Horv\'ath, {\it Several Problems of $n$-dimensional Discrete Geometry},
Ph.D Thesis, Steklov Math. Inst. Moscow, 1986.
\bibitem{kann88}R. Kannan and L. Lov\'asz, Covering minima and
lattice-point-free convex bodies, {\it Ann. Math.} {\bf 128} (1988),
577-602.
\bibitem{linh78}J. Linhart, Closest packings and closest coverings
by translates of a convex disc, {\it Studia Math. Hungar.} {\bf 13}
(1978), 157-162.
\bibitem{mahl47}K. Mahler, On the minimum determinant and the
circumscribed hexagons of a convex domain, {\it Proc. Kon. Ned.
Akad. Wet.} {\bf 50} (1947), 692-703.
\bibitem{rein34}K. Reinhardt, \"Uber die dichteste gitterf\"ormige
Lagerung kongruenter Bereiche in der Ebene und eine besondere Art
konvexer Kurven, {\it Abh. Math. Sem. Hamburg} {\bf 10} (1934),
216-230.
\bibitem{roge50}C.A. Rogers, A note on coverings and packings, {\it J. London Math. Soc.} {\bf
25} (1950), 327-331.
\bibitem{roge51}C.A. Rogers, The closest packing of convex two-dimensional domains,
{\it Acta Math.} {\bf 86} (1951), 309-321.
\bibitem{rysk70}S.S. Ry$\check{s}$kov, The polyhedron $\mu (m)$ and
certain extremal problems of the geometry of numbers (in Russian),
{\it Dokl. Akad. Nauk SSSR} {\bf 194} (1970), 514-517.
\bibitem{zong96}C. Zong, A few remarks on kissing numbers of a
convex body, {\it Anz. \"Osterreich. Akad. Wiss. Math.-Naturwiss.
KL.} {\bf 132} (1996), 11-15.
\bibitem{zong02}C. Zong, Simultaneous packing and covering in the
Euclidean plane, {\it Monatsh. Math.} {\bf 134} (2002), 247-255.
\bibitem{zong02'}C. Zong, From deep holes to free planes, {\it Bull. Amer.
Math. Soc.} {\bf 39} (2002), 533-555.
\bibitem{zong03}C. Zong, Simultaneous packing and covering in
three-dimensional Euclidean space, {\it J. London Math. Soc.} {\bf
67} (2003), 29-40.

\end{thebibliography}

\end{document}